\documentclass[11pt,openany]{article}

\usepackage{graphicx} 
\usepackage{subfig}
\usepackage{amssymb}
\usepackage{amsmath}
\usepackage{amsthm}
\usepackage{amsfonts}
\usepackage{mathrsfs}
\usepackage{makeidx}
\usepackage{multicol}

\input{cyracc.def}

\numberwithin{equation}{section}

\begin{document}
\title{\Large \bf  All meromorphic solutions of Fermat-type functional equations}
\author{ Feng L\"{u}\footnote{Corresponding author, Email: lvfeng18@gmail.com(F. L\"{u}).}\\
\small{College of Science, China University of Petroleum, Qingdao}\\
\small{Shandong, 266580, P.R. China}\\
}
\date{}
\maketitle
\vspace{3mm}

\begin{abstract}
In this paper, by making use of properties of elliptic functions, we describe meromorphic solutions of Fermat-type functional equations $f(z)^{n}+f(L(z))^{m}=1$ over the complex plane $\mathbb{C}$, where $L(z)$ is a nonconstant entire function, $m$ and $n$ are two positive integers. As applications, we also consider meromorphic solutions of Fermat-type difference and $q$-difference equations.
\end{abstract}



\pagenumbering{arabic}

\section{Introduction and main results}
\quad In 1637, Fermat \cite{Di} stated the conjecture (which is known as Fermat's last theorem) that the equation $x^m + y^m =1$ cannot have positive rational solutions if $m>2$. Since then, the equation has been a subject of intense and often heated discussions. In 1995, Wiles \cite{W2, W1} proved the profound conjecture.\\

In 1927, Montel in \cite{PM} initially considered the functional equations
\begin{equation}\label{1.1}
f^m(z)+g^m(z)=1,
\end{equation}
which can be regarded as the analogous of Fermat diophantine equations $x^m + y^m =1$ over function fields. He showed that all the entire solutions $f(z)$ and $g(z)$ of (\ref{1.1}) must be constant if $m\geq 3$, see also Jategaonka \cite{Ja}. The follow-up works were given by Baker in \cite{Ba} and Gross in \cite{FG}, respectively, they generalized Montel's result by proving that (\ref{1.1}) does not have nonconstant meromorphic solutions when $m>4$ and described nonconstant meromorphic solutions for $n=2,~3$. In 1970, Yang \cite{Y} considered the more general functional equations
\begin{equation}\label{1.2}
f^n(z)+g^m(z)=1,
\end{equation}
and derived that $f^m(z)+g^n(z)=1$ does not admit nonconstant entire solutions if $\frac{1}{n}+\frac{1}{m}<1$. Since then, (\ref{1.2}) has been studied in various settings, see \cite{G, Li1, Li2, Y}. For the convenience, some results can be stated as follows. (see e.g., \cite[Proposition 1]{Chen}, \cite{Ba, F.G, Yana}).\\

\noindent \textbf{Theorem A.} Suppose $f(z)$ and $g(z)$ are nonconstant meromorphic solutions of the functional equation (\ref{1.2}). \\

(i) If $n=m=3$, then $f(z)=\frac{\frac{1}{2}\left\{1+\frac{\wp^{\prime}(h(z))}{\sqrt{3}}\right\}}{\wp(h(z))}$ and $g(z)=\frac{\frac{\eta}{2}\left\{1-\frac{\wp^{\prime}(h(z))}{\sqrt{3}}\right\}}{\wp(h(z))}$ for nonconstant entire function $h$,  where $\eta^{3}=1$ and $\wp$ denotes the Weierstrass $\wp$-function satisfying $(\wp')^2\equiv 4\wp^3-1$ after appropriately choosing its periods. \\

(ii) If $n = 2$ and $m=3$, then $f(z) = i\wp'(h(z))$ and $g(z) = \eta\sqrt[3]{4}\wp(h(z)),$ where $h(z)$ is any nonconstant entire function. \\

(iii)If $m=2$ and $n=4$ , then $f(z)=sn'(h(z))$ and $g(z)= sn(h(z))$, where $h(z)$ is any nonconstant entire function and $sn$ is Jacobi elliptic function satisfying ${sn'}^2=1- sn^4$.\\

In 1989, Yanagihara \cite{Yana} considered the existence of meromorphic solutions of Fermat-type functional equation in another direction. In fact, Yanagihara obtained that $f^3(z)+f^3(z+c)=1$ does not admit nonconstant meromorphic functions of finite order. This result was also derived by L\"{u}-Han in \cite{LH} by making use of the difference analogue of the logarithmic derivative lemma of finite order meromorphic functions, which was established by Halburd and Korhonen \cite{HK}, Chiang and Feng \cite{CF}, independently. In 2014, this difference analogue was improved by Halburd, Korhonen and Tohge \cite{HAL} to meromorphic functions of hyper-order strictly less than 1. Here, the order and hyper-order of a meromorphic function $f$ are defined as
$$
\rho(f):=\limsup\limits_{r\to+\infty}\frac{\log T(r,f)}{\log r},~~\rho_{2}(f)={\limsup_{r\to +\infty}}\frac{\log\log T(r,f)}{\log r},
$$
where $T(r,f)$ denotes the Nevanlinna characteristic function of $f$. \\

With the help of these results, Korhonen and Zhang in \cite{KZ} generalize Yanagihara's theorem to meromorphic solutions of hyper-order strictly less than 1 as follows.\\

\noindent \textbf{Theorem B.} The functional equation
$$
f^3(z)+f^3(z+c)=1
$$
does not admit nonconstant meromorphic solutions of hyper-order strictly less than 1.\\

Recently, L\"{u} and his co-workers \cite{BL,LH, LG} also considered this kind of problems, and their results can be stated as follows.\\

\noindent \textbf{Theorem C.} Let $P(z)$ be a polynomial and $f(z)$ be a nonconstant meromorphic function of hyper-order strictly less than 1. If $f(z)$ is a solution of
\begin{equation}\label{1.3}
f^3(z)+f^3(z+c)=e^{P(z)},
\end{equation}
then, $P(z)=\alpha z+\beta$ is a linear function and $f(z)=de^{\frac{\alpha z+\beta}{3}}$, where $d(\neq 0)$ is constant and $d^3(1+e^{\alpha c})=1$. \\

In recent years, Fermat-type difference as well as differential-difference equations have been
studied extensively. As a result, successively a lot of investigations have done
by many scholars in this direction, see e.g., \cite{Liu3, Liu4}. We note that most of above results, including Theorems B and C, were obtained under the condition that the solutions is of hyper-order strictly less than 1, since the difference analogue of the logarithmic derivative lemma is needed in proofs of them. So, it is natural to ask what will happen if the hyper-order condition is omitted. We find that the conclusions of some previous theorems maybe invalid. In \cite{LH}, Han offered an example to show this point as follows.\\

\textbf{Example 1.} Let $c=\pi i$ and $\alpha, \beta$ be fixed constants satisfying $e^{\alpha c}=1$. Consider $f(z)=\frac{1}{2}\frac{1+\frac{\wp'(h(z))}{\sqrt{3}}}{\wp(h(z))}e^{\frac{\alpha z+\beta}{3}}$, where $h(z)=e^z$.  Then, a routine computation leads to $ f(z+c)=\frac{\eta}{2}\frac{1-\frac{\wp'(h(z))}{\sqrt{3}}}{\wp(h(z))}e^{\frac{\alpha z+\beta}{3}}$, where $\eta=e^{\frac{\alpha c}{3}}$.  Further, one has
 $$
 f^3(z)+f^3(z+c)=e^{\alpha z+\beta},
 $$
which implies that the above equation admits meromorphic solution $f(z)$. Obviously, $f(z)$ is not the form $de^{\frac{\alpha z+\beta}{3}}$ and the hyper-order of $f$ is $\rho_2(f)=1$.\\

Therefore, if the hyper-order condition is omitted, the Fermat-type difference equation may admit the other type of meromorphic solutions. In this paper, we pay attention to above question. Due to properties of elliptic functions, we describe the forms of all meromorphic solutions of some Fermat-type difference equations.\\

Before giving our main results, we introduce the Weierstrass $\wp$-function.  \\

The Weierstrass $\wp$-function is elliptic (also doubly periodic) function with periods $\omega_{1}$ and $\omega_{2}$ (${\sf Im}\frac{\omega_{1}}{\omega_{2}}\neq 0$) which is defined as
$$
\wp(z)=\wp\left(z ; \omega_{1}, \omega_{2}\right):=\frac{1}{z^{2}}+\sum_{\mu, \nu \in \mathbf{Z} ; \mu^{2}+\nu^{2} \neq 0}\left\{\frac{1}{\left(z+\mu \omega_{1}+\nu \omega_{2}\right)^{2}}-\frac{1}{\left(\mu \omega_{1}+\nu \omega_{2}\right)^{2}}\right\},
$$
and satisfies, after appropriately choosing $\omega_{1}$ and $\omega_{2}$,
$$
\left(\wp^{\prime}\right)^{2}=4 \wp^{3}-1.
$$
The period of $\wp$ span the lattice $\mathbb{L}=\omega_{1}\mathbb{Z}\oplus\omega_{2}\mathbb{Z}$. We also denote
$$
\mathbb{L}=\{\mu \omega_{1}+\nu \omega_{2}:~~ \mu,~\nu=0,\pm1, \pm2,...\}.
$$
Obviously, all the points in $\mathbb{L}$ are the poles and periods of $\wp$. Suppose that $D$ is the parallelogram with vertices at 0, $\omega_1$, $\omega_2$, $\omega_1+\omega_2$. Note that the order of $\wp$ is 2. Here, the order is the number of poles of $\wp$ or the number of zeros of $\wp-a$ ($a\in \mathbb{C}$) in the parallelogram. Together with $(\wp')^2=4 \wp^{3}-1$, one gets that $\wp$ has two distinct zeros in $D$, say $\theta_1$ and $\theta_2$. In view of $\wp(z)=0 \Leftrightarrow \wp'(z)=\pm i$, without loss of generality, throughout the paper, we assume that $\wp'(\theta_1)=-i$ and $\wp'(\theta_2)=i$. \\

Here, for two meromorphic functions $f,~g$ and two points $a,~b\in \mathbb{C}\cup\{\infty\}$, the notation $f(z)-a=0\Rightarrow g(z)-b=0$ means that all the zeros of $f-a$ are the zeros of $g(z)-b$. And the notation $f(z)-a=0\Leftrightarrow g(z)-b=0$ means that $f-a$ and $g-b$ have the same zeros.\\

In the present paper, more generally, we characterize meromorphic solutions of Fermat-type functional equations as follows.\\

\noindent \textbf{Proposition 1.} Suppose that $L(z)$ is a nonconstant entire function. Then, the nonconstant meromorphic function $f(z)$ is a solution of
\begin{equation}\label{1.002}
f^3(z)+f^3(L(z))=1
\end{equation}
if and only if $L(z)=qz+c$ is a linear function with $|q|=1$, and $f(z)=\frac{1}{2}\left\{1+\frac{\wp^{\prime}(h(z))}{\sqrt{3}}\right\} / \wp(h(z))$, where $h(z)$ is a nonconstant entire function satisfies one of the following equations.

(1). $h(L(z))=Ah(z)+\theta_1+\tau_1$ with $\tau_1\in \mathbb{L}$;

(2). $h(L(z))=Ah(z)+\theta_2+\tau_2$ with $\tau_2\in \mathbb{L}$;

(3). $h(L(z))=Ah(z)+\tau_3$ with $\tau_3\in \mathbb{L}$, where $A$ is a constant with $A^3=-1$ and $A\tau\in \mathbb{L}$ for any period of $\wp(z)$.\\

\noindent \textbf{Proposition 2.} Suppose that $L(z)$ is a nonconstant entire function, suppose that $g(z)$ is an entire function and $f(z)$ is a nonconstant meromorphic function, and suppose that $n(\geq 2),~m(\geq 3)$ are two integers such that $(n,m)\neq (3,3)$. Then, $f(z)$ is a solution of
\begin{equation}\label{L}
f^n(z)+f^m(L(z))=e^{g(z)}
\end{equation}
if and only if $L(z)=qz+c$ is a linear function, $f(z)=Ae^\frac{g(z)}{n}$ and
\begin{equation}\label{L0111}
g(L(z))=\frac{n}{m}g(z)+nLn \frac{B}{A},
\end{equation}
where $A,~B$ are two nonzero constants with $A^n+B^m=1$. In particularly, if $m=n$, then $|q|=1$.\\

It is pointed out that if $m\neq n$ in Proposition 2, the conclusion $|q|=1$ maybe invalid, as shown by the following example. \\

Suppose $m\neq n$ and $q^2=\frac{m}{n}$. Consider $g(z)=z^2+a$ with a constant $a$ satisfying $(1-\frac{m}{n})a=nLn \frac{B}{A}$. Then, a calculation yields $g(qz)=\frac{m}{n}g(z)+nLn \frac{B}{A}$ and $f(z)=Ae^\frac{g(z)}{n}$ is a solution of (\ref{L}). But $|q|\neq1$.\\

Let's turn back to the Fermat-type difference equations. Suppose $L(z)=z+c$ with $c\neq 0$. Then, by Proposition 1, we obtain the following result.\\

\noindent \textbf{Theorem 1.} The nonconstant meromorphic function $f(z)$ is a solution of
\begin{equation}\label{1.0002}
f^3(z)+f^3(z+c)=1
\end{equation}
if and only if $f(z)=\frac{1}{2}\left\{1+\frac{\wp^{\prime}(h(z))}{\sqrt{3}}\right\} / \wp(h(z))$, where $h(z)$ is a nonconstant entire function satisfies one of the following equations.

(1). $h(z+c)=Ah(z)+\theta_1+\tau_1$ with $\tau_1\in \mathbb{L}$;

(2). $h(z+c)=Ah(z)+\theta_2+\tau_2$ with $\tau_2\in \mathbb{L}$;

(3). $h(z+c)=Ah(z)+\tau_3$ with $\tau_3\in \mathbb{L}$,
where $A$ is a constant with $A^3=-1$ and $A\tau\in \mathbb{L}$ for any period of $\wp(z)$.\\

\textbf{Remark 1.} It is known that if $h(z)$ satisfies one of (1)-(3) in Theorem 1, then the order $\rho(h(z))\geq 1$ and hyper-order $\rho_2(f)\geq 1$. (The fact can be found in \cite{G}). So, if $\rho_2(f)<1$, then (\ref{1.002}) does not admit nonconstant meromorphic solutions. This is the conclusion of Theorem B. Below, we offer example to show that there exist entire function $h(z)$ with arbitrary order $\rho(h(z))(\geq 1)$ satisfying one of (1)-(3).\\

\textbf{Example 2.} In \cite[Theorem 3]{O}, Ozawa derived that there exists periodic entire function of arbitrary order $\sigma(\geq 1)$. So, there exists an entire function $g(z)$ such that $\rho(g(z))=\sigma$ and $g(z+c)=g(z)$. Set $h(z)=e^{az}g(z)+b$, where $a$, $b$ are constants with $e^{ac}=A$ and $(1-A)b=\theta_1+\tau_1$ or $\theta_2+\tau_2$ or $\tau_3$ with $A^3=-1$. Then, $\rho(h(z))=\rho(g(z))=\sigma$ and a calculation yields that
$$
h(z+c)=Ah(z)+(1-A)b.
$$
Thus, $h(z)$ satisfies one of (1)-(3) in Theorem 1.\\

By Theorem 1 and Theorem C, we can get the following result. \\

\noindent \textbf{Theorem 2.} Suppose that $f(z)$ is a nonconstant meromorphic function and $P(z)$ is a polynomial. Then $f(z)$ is a solution of the functional equation
\begin{equation}\label{Y11}
f^3(z)+f^3(z+c)=e^{P(z)}
\end{equation}
if and only if $P(z)=\alpha z+\beta$ with two constants $\alpha,~\beta$ and $f(z)$ satisfies one of the following cases.

(a). $f(z)=de^{\frac{\alpha z+\beta}{3}}$, where $d$ is a nonzero constant and $d^3(1+e^{\alpha c})=1$;

(b). $f(z)=e^{\frac{\alpha z+\beta}{3}}\frac{1}{2}\left\{1+\frac{\wp^{\prime}(h(z))}{\sqrt{3}}\right\} / \wp(h(z))$, where $e^{\alpha c}=1$ and $h(z)$ is a nonconstant entire function satisfies one of (1)-(3) in Theorem 1.\\

With Proposition 2, we derive a theorem as follows.\\

\noindent \textbf{Theorem 3.} Suppose that $g(z)$ is an entire function and $f(z)$ is a nonconstant meromorphic function. Suppose that $n(\geq 2),~m(\geq 3)$ are two integers such that $(n,m)\neq (3,3)$. Then, $f(z)$ is a solution of
\begin{equation}\label{L1}
f^n(z)+f^m(z+c)=e^{g(z)}
\end{equation}
if and only if $f(z)=Ae^\frac{g(z)}{n}$ and $g(z+c)=\frac{n}{m}g(z)+nLn \frac{B}{A}$, where $A,~B$ are two nonzero constants with $A^n+B^m=1$. In particularly, if the order of $g(z)$ is less than 1, then (\ref{L1}) admits nonconstant meromorphic solutions if and only if $n=m$ and $g(z)=\alpha z+\beta$, $f(z)=Ae^\frac{\alpha z+\beta}{n}$ with constant $\alpha(\neq0)$ such that $A^n(1+e^{\alpha c})=1$.\\

\textbf{Remark 2.} We point out that equation (\ref{L1}) may have nonconstant meromorphic functions if $\rho(g(z))\geq 1$, as shown by the following examples.\\

\textbf{Example 3.} Suppose that $m\neq n$, $g(z)=e^{\alpha z+\beta}+a$ with constants $\alpha(\neq0),~\beta,~a$. Consider $f(z)=Ae^{\frac{g(z)}{n}}$ with $e^{\alpha c}=\frac{n}{m}$ and $(1-\frac{n}{m})a=n Ln \frac{B}{A}$, where $A,~B$ are nonzero constants with $A^n+B^m=1$. By a calculation, one gets
$g(z+c)=\frac{n g(z)}{m}+n Ln \frac{B}{A}$ and
$$
f^n(z)+f^m(z+c)=e^{g(z)}.
$$
Obviously, $\rho(g(z))=1$.\\

\textbf{Example 4.} Suppose that $m=n$, and $h$ is any nonconstant entire periodic function with period $c$. Let $g(z)= h(z)+az$ with $ac=nLn\frac{B}{A}$, where $A,~B$ are nonzero constants with $A^n+B^n=1$. Then, we have $g(z+c)=g(z)+nLn\frac{B}{A}$, and $f(z)=Ae^{\frac{g(z)}{n}}$ is a nonconstant meromorphic solution of
$$
f^n(z)+f^n(z+c)=e^{g(z)}.
$$
Obviously, $\rho(g(z))\geq 1$.\\

\textbf{Remark 3.} From Theorem 3, we see that the equation (\ref{L1}) do not admit meromorphic solutions with poles. It is pointed out that the same argument in Theorems 3 can deal with the Fermat-type difference equations $f^m(z)+f^n(z+c)=e^{g(z)}$ with $n(\geq 2),~m(\geq 3)$ and $(m,n)\neq (3,3)$, we omit the details here. \\

The following corollary follows immediately from Theorem 3.\\

\noindent \textbf{Corollary 1.} Suppose that $P(z)$ is a polynomial and $f(z)$ is a nonconstant meromorphic function. Suppose that $n(\geq 2),~m(\geq 2)$ are two integers with $m+n\geq5$ and $(n,m)\neq (3,3)$. Then, $f(z)$ is a solution of
\begin{equation}\label{LL11}
f^n(z)+f^m(z+c)=e^{P(z)}
\end{equation}
if and only $n=m$, $f(z)=Ae^\frac{P(z)}{n}$ and $P=\alpha z+\beta$ with two constants $\alpha(\neq0),~\beta$. \\

From Theorem 3 and Corollary 1, one can easily get the following result.\\

\noindent \textbf{Corollary 2.} Suppose that $f(z)$ is a nonconstant meromorphic function, and $n(\geq 2),~m(\geq 2)$ are two integers with $m+n\geq5$. Then, $f(z)$ is a solution of
\begin{equation}\label{LL21}
f^n(z)+f^m(z+c)=1
\end{equation}
if and only if $n=m=3$. \\

Finally, as applications of Propositions 1 and 2, we consider meromorphic solutions of Fermat-type $q$-difference functional equations. As early as 1952, Valiron in \cite{V} showed that the non-autonomous Schr\"{o}der $q$-difference equation
\begin{equation}\label{M1}
f(qz)=R(z,f(z))
\end{equation}
where $R(z, f(z))$ is rational in both arguments, admits a one parameter family of
meromorphic solutions, provided that $q\in \mathbb{C}$ is chosen appropriately. Later, Gundersen et al \cite{Gun} proved that if $|q|>1$ and the $q$-difference equation (\ref{M1}) admits a meromorphic solution of order zero, then (\ref{M1}) reduces to
a $q$-difference Riccati equation, i.e. $\deg_f R = 1$. Some scholars, such as Bergweiler-Hayman \cite{Ber1}, Bergweiler-Ishizaki-
Yanagihara \cite{Ber2}, Eremenko-Sodin \cite{Ere}, Ishizaki-Yanagihara \cite{Is} also made contributions to meromorphic solutions of $q$-difference functional equations. In 2007, Barneet-Halburd-Korhonen-Morgan in \cite{Bar} derived the $q$-difference analogues of some well-known results in  Nevanlinna theory, including the lemma of the logarithmic derivative, the Clunie's lemma and the second main theorem. With the help of results, meromorphic solutions to $q$-difference functional equations were further studied, see \cite{KW, Liu5}.\\

Observe that in Propositions 1 and 2, $L(z)$ reduce to a linear function, say $qz+c$ with $q\neq 0$. Therefore, we have described  meromorphic solutions of Fermat-type $q$-difference functional equations $f^3(z)+f^3(qz+c)=1$ and $f^n(z)+f^n(qz+c)=e^{g(z)}$, where $n,~m$ satisfy some condition and $g(z)$ is an entire function. \\

\textbf{Remark 4.} From the above theorems, we see that the case $(1,n),~(m,1),~(2,2)$ are left. Unfortunately, we cannot deal with these cases and thus leave them for further study.\\

For the proofs, we will assume that the reader is familiar with basic elements in Nevanlinna theory of
meromorphic functions in $\mathbb{C}$ (see e.g. \cite{HAL, Hay, LI, YY}), such as
the {\it first} and {\it second} main theorems, the {\it
characteristic function} $T(r,f)$, the {\it proximity function}
$m(r,f)$, the {\it counting function} $N(r,f)$ and the {\it reduced counting function} $\overline{N}(r,f)$. We also need the following notation.\\

The lower order of a meromorphic function $f$ is defined as
$$
\mu(f):=\liminf\limits_{r\to+\infty}\frac{\log T(r,f)}{\log r}.
$$

\section{Proofs of main results}

In this section, we firstly give the proof of Proposition 1.

\begin{proof}[Proof of Proposition 1]Firstly, we will prove the necessity. Assume that $h(z)$ satisfies one of (1)-(3), we will prove that $f(z)$ is a meromorphic solution of (\ref{1.002}). Suppose that $h(z)$ satisfies (1). Then $h(L(z))=Ah(z)+\theta_1+\tau_1$ with $\tau_1\in \mathbb{L}$. The assumption that $A\tau\in \mathbb{L}$ for any period of $\wp(z)$ defines a mapping as
$$
A:\mathbb{L}\rightarrow \mathbb{L},~~A:\tau\rightarrow A\tau.
$$
We claim that the mapping $A:$ is a bijection. \\

Observe that $A=-1$, or $e^{\frac{\pi}{3}i}$, or $e^{-\frac{\pi}{3}i}$. If $A=-1$, clearly, the claim is valid. Suppose that $A=e^{\frac{\pi}{3}i}$. Among all the periods of $\wp(\omega)$, we list the periods which may take the smallest modulus as follows: $\omega_1$, $\omega_1+\omega_2$, $\omega_2$, $\omega_2-\omega_1$, $-\omega_1$, $-\omega_2-\omega_1$, $-\omega_2$ and $\omega_1-\omega_2$. Then, in view of that $A\tau=\tau e^{\frac{\pi}{3}i}$ is also a period of $\wp$, below, we have two possible cases (see Figures 1-2 below).\\

\begin{figure}[htbp]
\begin{minipage}[t]{0.5\linewidth}
\centering
\includegraphics[scale=0.7]{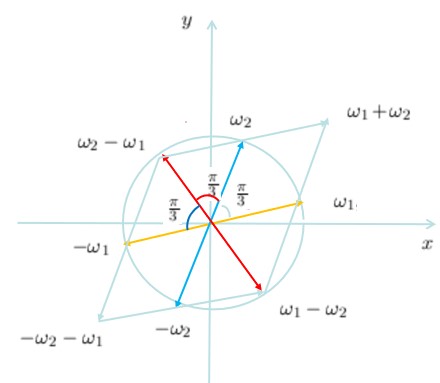}
\caption{(i)}
\end{minipage}%
\begin{minipage}[t]{0.5\linewidth}
\centering
\includegraphics[scale=0.7]{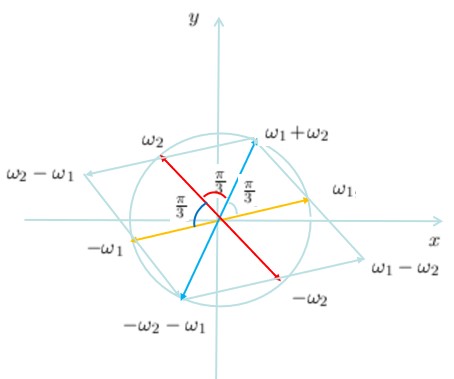}
\caption{(ii)}
\end{minipage}
\end{figure}

\textbf{(i).} $A\omega_1=\omega_2$ and $A\omega_2=\omega_2-\omega_1$ and $A(\omega_2-\omega_1)=-\omega_1$.\\

\textbf{(ii).} $A\omega_1=\omega_1+\omega_2$ and $A(\omega_1+\omega_2)=\omega_2$ and $A\omega_2=-\omega_1$.\\

Either (i) or (ii) holds, a simple analysis yields the claim is right. The same argument implies that the claim holds for $A=e^{-\frac{\pi}{3}i}$. Therefore, the claim is proved.\\

Observe that $A^3=-1$, the claim and the definition of $\wp(z)$ as follows
$$
\wp(z)=\wp\left(z ; \omega_{1}, \omega_{2}\right):=\frac{1}{z^{2}}+\sum_{\mu, \nu \in \mathbf{Z} ; \mu^{2}+\nu^{2} \neq 0}\left\{\frac{1}{\left(z+\mu \omega_{1}+\nu \omega_{2}\right)^{2}}-\frac{1}{\left(\mu \omega_{1}+\nu \omega_{2}\right)^{2}}\right\}.
$$
A calculation yields that
\begin{equation}\label{UV1}
\begin{aligned}
\wp(Az)&=\frac{1}{(Az)^{2}}+\sum_{\mu, \nu \in \mathbf{Z} ; \mu^{2}+\nu^{2} \neq 0}\left\{\frac{1}{\left(Az+\mu \omega_{1}+\nu \omega_{2}\right)^{2}}-\frac{1}{\left(\mu \omega_{1}+\nu \omega_{2}\right)^{2}}\right\}\\
&=\frac{1}{A^2}\left\{\frac{1}{z^{2}}+\sum_{\mu, \nu \in \mathbf{Z} ; \mu^{2}+\nu^{2} \neq 0}\left\{\frac{1}{\{z+\frac{1}{A}(\mu \omega_{1}+\nu \omega_{2})\}^{2}}-\frac{1}{\{\frac{1}{A}(\mu \omega_{1}+\nu \omega_{2})\}^{2}}\right\}\right\}\\
&=-A\left\{\frac{1}{z^{2}}+\sum_{\mu, \nu \in \mathbf{Z} ; \mu^{2}+\nu^{2} \neq 0}\left\{\frac{1}{\{z-A^2(\mu \omega_{1}+\nu \omega_{2})\}^{2}}-\frac{1}{\{-A^2(\mu \omega_{1}+\nu \omega_{2})\}^{2}}\right\}\right\}\\
&=-A\left\{\frac{1}{z^{2}}+\sum_{\mu, \nu \in \mathbf{Z} ; \mu^{2}+\nu^{2} \neq 0}\left\{\frac{1}{\left(z+\mu \omega_{1}+\nu \omega_{2}\right)^{2}}-\frac{1}{\left(\mu \omega_{1}+\nu \omega_{2}\right)^{2}}\right\}\right\}\\
&=-A\wp(z),
\end{aligned}
\end{equation}
which implies that $\wp'(Az)=-\wp'(z)$. It is known that
$$
\wp(w+c)=\frac{1}{4}[\frac{\wp'(w)-\wp'(c)}{\wp(w)-\wp(c)}]^2-\wp(w)-\wp(c),
$$
for any $w$ and $c$. Thus,
\begin{equation}\label{X2.13}
\begin{aligned}
\wp(Aw+\theta_1)&=\frac{1}{4}[\frac{\wp'(Aw)-\wp'(\theta_1)}{\wp(Aw)-\wp(\theta_1)}]^2-\wp(Aw)-\wp(\theta_1)=\frac{1}{4}[\frac{\wp'(Aw)+i}{\wp(Aw)}]^2-\wp(Aw)\\
&=\frac{1}{4}\frac{(\wp'(Aw)+i)^2\wp(Aw)}{\wp(Aw)^3}-\wp(Aw)=\frac{(\wp'(Aw)+i)^2\wp(Aw)}{\wp'(Aw)^2+1}-\wp(Aw)\\
&=\frac{(\wp'(Aw)+i)^2\wp(Aw)}{(\wp'(Aw)+i)(\wp'(Aw)-i)}-\wp(Aw)=\frac{2i\wp(Aw)}{\wp'(Aw)-i}\\
&=\frac{-A2i\wp(w)}{-\wp'(w)-i}=\frac{A2i\wp(w)}{\wp'(w)+i}.\\
\end{aligned}
\end{equation}
\begin{equation}\label{X2.14}
\begin{aligned}
\wp'(Aw+\theta_1)&=(\frac{2i\wp(w)}{\wp'(w)+i})'=2i\frac{\wp'(w)(\wp'(w)+i)-\wp(w)\wp''(w)}{(\wp'(w)+i)^2}\\
&=2i\frac{\wp'(w)^2+i\wp'(w)-\wp(w)\wp''(w)}{(\wp'(w)+i)^2}=2i\frac{\wp'(w)^2+i\wp'(w)-\frac{3}{2}(\wp'(w)^2+1)}{(\wp'(w)+i)^2}\\
&=-i\frac{\wp'(w)-3i}{\wp'(w)+i}.
\end{aligned}
\end{equation}
Set $h(z)=\omega$. Then $h(L(z))=Ah(z)+\theta_1+\tau_1=Aw+\theta_1+\tau_1$. Note that both $\wp$ and $\wp'$ are elliptic functions with periods $\omega_{1}$ and $\omega_{2}$, we have
$$
\begin{aligned}
&\wp(h(L(z)))=\wp(Aw+\theta_1+\tau_1)=\wp(Aw+\theta_1),\\
&\wp'(h(L(z)))=\wp'(Aw+\theta_1+\tau_1)=\wp'(Aw+\theta_1).
\end{aligned}$$
Observe that $f(z)=\frac{1}{2}\frac{1+\frac{\wp'(h(z))}{\sqrt{3}}}{\wp(h(z))}$. So
\begin{equation}\label{X2.15}
\begin{aligned}
f(L(z))&=\frac{1}{2}\frac{1+\frac{\wp'(h(L(z)))}{\sqrt{3}}}{\wp(h(L(z)))}=\frac{1}{2}\frac{1+\frac{\wp'(Aw+\theta_1)}{\sqrt{3}}}{\wp(Aw+\theta_1)}=\frac{1}{2}\frac{1-\frac{1}{\sqrt{3}}i\frac{\wp'(w)-3i}{\wp'(w)+i}}{\frac{A2i\wp(w)}{\wp'(w)+i}}\\
&=\frac{1}{2}\frac{(\wp'(w)+i)-\frac{1}{\sqrt{3}}i(\wp'(w)-3i)}{A2i\wp(w)}
=\frac{1}{2}\frac{-1-\sqrt{3}i}{2}\frac{1}{-A}\frac{1-\frac{\wp'(w)}{\sqrt{3}}}{\wp(w)}=\frac{1}{-A}\frac{\eta}{2}\frac{1-\frac{\wp'(w)}{\sqrt{3}}}{\wp(w)},\\
\end{aligned}
\end{equation}
where $\eta=\frac{-1-\sqrt{3}i}{2}$ and $\eta^3=1$. Thus, by (i) of Theorem A, one can easily derive that $f(z)^3+f(L(z))^3=1$ and $f(z)$ is a solution of (\ref{1.002}). If $h(z)$ satisfies (2) or (3), then the same argument yields that $f(z)$ is also a solution of (\ref{1.002}). \\

Conversely, we will prove the sufficiency. For convenience, we below denote $L(z)$ by $\overline{z}$, which means that $\overline{z}=L(z)$ and a point $\overline{a}=L(a)$.\\

Suppose that $f(z)$ is a nonconstant meromorphic solution of (\ref{1.002}). Then, Via (i) of Theorem A, one has
\begin{equation}\label{2.1}
f(z)=\frac{1}{2}\frac{1+\frac{\wp'(h(z))}{\sqrt{3}}}{\wp(h(z))} \qquad and \qquad f(\overline{z})=\frac{\eta}{2}\frac{1-\frac{\wp'(h(z))}{\sqrt{3}}}{\wp(h(z))},
\end{equation}
where $h$ is a nonconstant entire function over $\mathbb{C}$. Next, we will prove that $h(z)$ satisfies one of (1)-(3). \\

Obviously, $T(r,f(z))=T(r,f(\overline{z}))+S(r,f)$. Rewrite the form of $f(z)$ as
$$
2f(z)\wp(h(z))-1=\frac{\wp'(h(z))}{\sqrt{3}}.
$$
Then,
\begin{equation}\label{X2.0001}
[2f(z)\wp(h(z))-1]^2=(\frac{\wp'(h(z))}{\sqrt{3}})^2=\frac{1}{3}[4\wp(h(z))^3-1],
\end{equation}
We rewrite (\ref{X2.0001}) as
\begin{equation}\label{X2.1}
f(z)[f(z)\wp(h(z))-1]=\frac{1}{3}\frac{[\wp(h(z))^3-1]}{\wp(h(z))},
\end{equation}
which implies that
$$
\begin{aligned}
3T(r,\wp(h(z)))&=T(r, \frac{1}{3}\frac{[\wp(h(z))^3-1]}{\wp(h(z))})+O(1)=T(r,f(z)[f(z)\wp(h(z))-1])+O(1)\\
&\leq T(r,f(z))+T(r,f(z)\wp(h(z)))+O(1)\\
&\leq 2T(r,f(z))+T(r,\wp(h(z)))+O(1).
\end{aligned}
$$
Thus, we derive that
\begin{equation}\label{X2.2}
T(r,\wp(h(z)))\leq T(r,f(z))+O(1).
\end{equation}
The form of $f(z)$ yields that $f(\overline{z})=\frac{1}{2}\frac{1+\frac{\wp'(h(\overline{z}))}{\sqrt{3}}}{\wp(h(\overline{z}))}$. Then, the same argument leads to
\begin{equation}\label{X2.21}
T(r,\wp(h(\overline{z})))\leq T(r,f(\overline{z}))+O(1).
\end{equation}
Rewrite (\ref{1.002}) as $f(\overline{z})^3=-[f(z)^3-1]=-(f(z)-1)(f(z)-\varsigma)(f(z)-\varsigma^{2}), ~(\varsigma\neq1,~\varsigma^3=1)$, which implies that
the zeros of $f(z)-1$, $f(z)-\varsigma$ and $f(z)-\varsigma^{2}$ are of multiplicities at least 3. Applying Nevanlinna's first and second theorems to the function $f(z)$ yields that
$$
\begin{aligned}
&2T(r,f(z))\leq\sum_{m=0}^{2}\overline{N}(r,\frac{1}{f(z)-\varsigma^{m}})+\overline{N}(r,f(z))+S(r,f(z)),\\
&\leq\frac{1}{3}\sum_{m=0}^{2} N(r,\frac{1}{f(z)-\varsigma^{m}})+N(r,f(z))+S(r,f(z))\\
&\leq2T(r,f(z))+S(r,f(z)).
\end{aligned}
$$
Therefore
\begin{equation}\label{X2.22}
T(r,f(z))=\overline{N}(r,f(z))+S(r,f(z))=N(r,f(z))+S(r,f(z)).
\end{equation}
Further, in view of that $f(z)=\frac{1}{2}\frac{1+\frac{\wp'(h(z))}{\sqrt{3}}}{\wp(h(z))}$ and $\wp$ only has multiple poles, one derives that
\begin{equation}\label{X2.3}
\begin{aligned}
T(r,f(z))&=\overline{N}(r,f(z))+S(r,f(z))=\overline{N}(r,\frac{1}{2}\frac{1+\frac{\wp'(h(z))}{\sqrt{3}}}{\wp(h(z))})+S(r,f(z))\\
&\leq \overline{N}(r,\frac{1}{\wp(h(z))})+\overline{N}(r,\wp(h(z)))+S(r,f(z))\\
&\leq \overline{N}(r,\frac{1}{\wp(h(z))})+\frac{1}{2}N(r,\wp(h(z)))+S(r,f(z))\\
&\leq \frac{3}{2}T(r,\wp(h(z)))+S(r,f(z)).
\end{aligned}
\end{equation}

The equation $f(z)^3+f(\overline{z})^3=1$ yields that $\overline{N}(r,f(\overline{z}))=\overline{N}(r,f(z))$. Then, by (\ref{X2.22}), one has
\begin{equation}\label{X2.4}
\begin{aligned}
T(r,f(\overline{z}))&=T(r,f(z))+S(r,f(z))=\overline{N}(r,f(z))+S(r,f(z))\\
&=\overline{N}(r,f(\overline{z}))+S(r,f(z)).
\end{aligned}
\end{equation}

Then, the same argument as in (\ref{X2.3}) yields that
\begin{equation}\label{X2.23}
T(r,f(\overline{z}))\leq \frac{3}{2}T(r,\wp(h(\overline{z})))+S(r,f(z)).
\end{equation}

All the above discussion yields that
$$S(r,f(z))=S(r,f(\overline{z}))=S(r,\wp(h(z)))=S(r,\wp(h(\overline{z}))).$$ For simplicity, we write $$S(r)=S(r,f(z))=S(r,f(\overline{z}))=S(r,\wp(h(z)))=S(r,\wp(h(\overline{z}))).$$

By (\ref{2.1}) and a routine computation, we get that
\begin{equation}\label{2.2}
\eta(1-\frac{\wp'(h(z))}{\sqrt{3}})\wp(h(\overline{z}))=(1+\frac{\wp'(h(\overline{z}))}{\sqrt{3}})\wp(h(z)).
\end{equation}

We employ the method in \cite{LH, FL, WL} to prove this theorem. For a set $G$, define the function $N(r, G)$ as
$$
N(r, G)=\int_{0}^r\frac{n(t,G)-n(0,G)}{t}dt+n(0,G)\log r,
$$
where the notation $n(r, G)$ is the number of points in $G\cap \{|z|<r\}$, ignoring multiplicities. We define a set $S_1$ as
$$
S_1=\{z|\wp(h(z))=0~~and~~\wp(h(\overline{z}))=0\}.
$$
Arrange $S_1$ as $S_1=\{a_{s}\}_{s=1}^{\infty}$ and $a_{s}\rightarrow\infty$ as $s\rightarrow\infty$. Notice, when $\wp(h(a_s))=0$, then $[\wp'(h(a_s))]^{2}=-1$. Further, when $\wp(h(\overline{a_s}))=0$, then $[\wp'(h(\overline{a_s}))]^{2}=-1$.  Differentiate (\ref{2.2}) and apply substitution to observe that
$$
\begin{aligned}
\begin{aligned}
&\quad\eta (1-\frac{\wp'(h(a_s))}{\sqrt{3}})\wp'(h(\overline{a_s}))[h(\overline{z})]'(a_s)\\
&=(1+\frac{\wp'(h(\overline{a_s}))}{\sqrt{3}})\wp'(h(a_s)))h'(a_s).
\end{aligned}
\end{aligned}
$$
From which we have that one and the only one of the following situations occurs
$$
\left\{
\begin{aligned}
	&g_1(a_s)=\eta(1-i\frac{\sqrt{3}}{3})[h(\overline{z})]'(a_s)-(1+i\frac{\sqrt{3}}{3})h'(a_s)=0, \\
	&g_2(a_s)=\eta [h(\overline{z})]'(a_s)+h'(a_s)=0,\\
	&g_3(a_s)=\eta(1+i\frac{\sqrt{3}}{3})[h(\overline{z})]'(a_s)-(1-i\frac{\sqrt{3}}{3})h'(a_s)=0.
\end{aligned}
\right.
$$
Below, we consider two cases.\\

\textbf{Case 1.} $g_i(z)\not\equiv 0$ for any $i=1,~2,~3$.\\

Here, we employ a result of Clunie \cite{Clu}, which can be stated as follows.\\

\emph{Lemma 1. Let $f$ be a nonconstant entire function and let $g$ be a transcendental meromorphic function
in the complex plane, then $T(r, f) = S(r, g(f)))$ as $r\rightarrow\infty $.}\\

Note that $\wp$ is transcendental. By Lemma 1, we have that $T(r,h(z))=S(r, \wp(h(z)))=S(r)$ and $T(r, h(\overline{z}))=S(r, \wp(h(\overline{z})))=S(r)$. We also have
$$T(r,h'(z))=O(T(r, h(z))=S(r),~~T(r,[h(\overline{z})]')=O(T(r, h(\overline{z})))=S(r).$$
 So, $T(r,g_i)=S(r)$. Further, we have
\begin{equation}\label{3.4}
\begin{aligned}
N(r, S_1)&\leq \sum_i N(r, \frac{1}{g_i})\leq \sum_i T(r, \frac{1}{g_i})\\
&=\sum_i T(r, g_i)+O(1)=S(r).
\end{aligned}
\end{equation}
Assume that $E$ is the set of all zeros of $\wp(h)$, that is $E=\{z|\wp(h(z))=0\}$. Obviously, $S_1\subseteq E$. Put $S_2=E\backslash S_1$. For any $b\in S_2$, we see that $\wp(h(b))=0$ and $\wp(h(\overline{b}))\neq 0$. In view of $[\wp'(h(b))]^{2}=-1$ and (\ref{2.2}), one has $\wp(h(\overline{b}))=\infty$. All the above discussions yields that
\begin{equation}\label{3.5}
\begin{aligned}
\overline{N}(r,\frac{1}{\wp(h(z))})&=\overline{N}(r,E)=N(r,S_1)+\overline{N}(r,S_2)\\
&\leq\overline{N}(r,\wp(h(\overline{z})))+S(r).
\end{aligned}
\end{equation}
Suppose that $z_0$ is a zero of $\wp(h(z))$ with multiplicity $p$. Then $z_0$ is a zero of $(\wp(h(z)))'=\wp'(h(z))h'(z)$ with multiplicity $p-1$. The fact $\wp'(h(z_0))=\pm i$ yields that $z_0$ is a zero of $h'$ with multiplicity $p-1$. Thus,
\begin{equation}\label{X2.5}
N(r,\frac{1}{\wp(h(z))})=\overline{N}(r,\frac{1}{\wp(h(z))})+N(r,\frac{1}{h'(z)})=\overline{N}(r,\frac{1}{\wp(h(z))})+S(r).
\end{equation}
The equation (\ref{X2.22}) yields that $m(r,f(z))=S(r)$, since $T(r,f(z))=m(r,f(z))+N(r,f(z))$. Rewrite the form of $f(z)$ as
$$
\frac{1}{\wp(h(z))}=2f(z)-\frac{\frac{\wp'(h(z))}{\sqrt{3}}}{\wp(h(z))}=2f(z)-\frac{\frac{\wp'(h(z))h'(z)}{\sqrt{3}}}{\wp(h(z))}\frac{1}{h'(z)}.
$$
Then, applying the logarithmic derivative lemma, we have
\begin{equation}\label{X2.6}
\begin{aligned}
m(r, \frac{1}{\wp(h(z))})&\leq m(r, f(z))+m (r, \frac{(\wp(h(z)))'}{\wp(h(z))})+m(r, \frac{1}{h'(z)})+O(1)\\
&\leq S(r)+S(r, \wp(h(z)))+S(r)=S(r).
\end{aligned}
\end{equation}
Combining (\ref{X2.5}) and (\ref{X2.6}) yields that
\begin{equation}\label{X2.7}
\begin{aligned}
T(r, \wp(h(z)))&=T(r,\frac{1}{\wp(h(z))})+O(1)=N(r,\frac{1}{\wp(h(z))})+m(r, \frac{1}{\wp(h(z))})+O(1)\\
&=\overline{N}(r,\frac{1}{\wp(h(z))})+m(r, \frac{1}{\wp(h(z))})+O(1)\\
&=\overline{N}(r,\frac{1}{\wp(h(z))})+S(r).
\end{aligned}
\end{equation}
We also know that $[\wp'(h(\overline{z}))]^2=4\wp(h(\overline{z}))^3-1$. Then,
$$[\wp(h(\overline{z}))]'^2=[h(\overline{z})]'^2\wp'(h(\overline{z}))^2=4[h(\overline{z})]'^2\wp(h(\overline{z}))^3-[h(\overline{z})]'^2.$$
Rewrite it as
$$
\wp(h(\overline{z}))^2 \wp(h(\overline{z}))=\wp(h(\overline{z}))^3=\frac{1}{4[h(\overline{z})]'^2} [\wp(h(\overline{z}))]'^2-\frac{1}{4}.
$$
By $T(r, \frac{1}{4[h(\overline{z})]'^2})=S(r)$ and Clunie's lemma (see \cite{Clu1}), we obtain that
$m(r,\wp(h(\overline{z})))=S(r)$ and
\begin{equation}\label{X2.8}
T(r,\wp(h(\overline{z})))=m(r,\wp(h(\overline{z})))+N(r,\wp(h(\overline{z})))=N(r,\wp(h(\overline{z})))+S(r).
\end{equation}
Note that $\wp(h(\overline{z}))$ only has multiple poles. So, $\overline{N}(r,\wp(h(\overline{z})))\leq \frac{1}{2}N(r,\wp(h(\overline{z})))$. Further, combining (\ref{X2.3}), (\ref{X2.23}), (\ref{3.5}), (\ref{X2.7}) and (\ref{X2.8}), we have that
\begin{equation}\label{X2.9}
\begin{aligned}
\frac{2}{3}T(r,f(z))&\leq T(r, \wp(h(z)))+S(r)=\overline{N}(r,\frac{1}{\wp(h(z))})+S(r)\leq\overline{N}(r,\wp(h(\overline{z})))+S(r)\\
&\leq \frac{1}{2}N(r,\wp(h(\overline{z})))+S(r)\leq \frac{1}{2}T(r,\wp(h(\overline{z})))+S(r)\\
&\leq \frac{1}{2}T(r,f(\overline{z}))+S(r)=\frac{1}{2}T(r,f(z))+S(r),
\end{aligned}
\end{equation}
which is a contradiction. Thus, the case cannot occur.\\

Case 2. $g_i(z)\equiv 0$, for some $i\in \{1,2, 3\}$.\\

Firstly, we assume that $g_1(z)\equiv 0$. Then, $[h(\overline{z})]'=Ah'(z)$, where $A=\frac{(1+i\frac{\sqrt{3}}{3})}{\eta(1-i\frac{\sqrt{3}}{3})}$. Integrating this equation yields $h(\overline{z})=Ah(z)+B$, where $B$ is a fixed constant.\\

We know that $\wp(h(z))$ has infinitely many poles. Suppose that $\wp(h(b_0))=\infty$. The equation (\ref{2.2}) yields that $\wp(h(\overline{b_0}))=0$ or $\wp(h(\overline{b_0}))=\infty$. Assume $\wp(h(\overline{b_0}))=\infty$. We rewrite (\ref{2.2}) as
\begin{equation}\label{X2.10}
\frac{1+\frac{\wp'(h(\overline{z}))}{\sqrt{3}}}{\wp(h(\overline{z}))}[h(\overline{z})]' =\frac{\eta(1-\frac{\wp'(h(z))}{\sqrt{3}})}{\wp(h(z))}[h(\overline{z})]' =\frac{\eta(1-\frac{\wp'(h(z))}{\sqrt{3}})}{\wp(h(z))}Ah'(z).
\end{equation}
Further, we rewrite (\ref{X2.10}) as
\begin{equation}\label{X2.11}
\begin{aligned}
\frac{[h(\overline{z})]' }{\wp(h(\overline{z}))}+\frac{1}{\sqrt{3}}\frac{[\wp(h(\overline{z}))]'}{\wp(h(\overline{z}))}=A\eta [\frac{h'(z)}{\wp(h(z))}-\frac{1}{\sqrt{3}} \frac{[\wp(h(z))]'}{\wp(h(z))}].\\
\end{aligned}
\end{equation}
Note that $\wp(h(z))$ only has multiple poles. Then, assume that $b_0$ is a pole of $\wp(h(z))$ and $\wp(h(\overline{z}))$ with multiplicities $m(\geq 2)$ and $k(\geq 2)$, respectively. By taking the residue of both sides of (\ref{X2.11}) at $b_0$, we have
$$
-k\frac{1}{\sqrt{3}}=A\eta\frac{m}{\sqrt{3}},
$$
which implies that $A\eta=-\frac{k}{m}$. It contradicts with $A\eta=\frac{1+i\frac{\sqrt{3}}{3}}{(1-i\frac{\sqrt{3}}{3})}=\frac{\sqrt{3}+i}{(\sqrt{3}-i)}$.\\

So, $\wp(h(\overline{b_0}))=0$ and $\wp'(h(\overline{b_0}))=\pm i$. Suppose that $\wp'(h(\overline{b_0}))=i$. Assume that $b_0$ is a zero of $\wp(h(\overline{z}))$ with multiplicity $s$. We rewrite (\ref{X2.11}) as
\begin{equation}\label{X2.12}
\begin{aligned}
\frac{1}{\wp'(h(\overline{z}))}\frac{[\wp(h(\overline{z}))]'}{\wp(h(\overline{z}))}+\frac{1}{\sqrt{3}}\frac{[\wp(h(\overline{z}))]'}{\wp(h(\overline{z}))}
&=\frac{[h(\overline{z})]'}{\wp(h(\overline{z}))}+\frac{1}{\sqrt{3}}\frac{[\wp(h(\overline{z}))]'}{\wp(h(\overline{z}))}\\
&=A\eta [\frac{h'(z)}{\wp(h(z))}-\frac{1}{\sqrt{3}} \frac{[\wp(h(z))]'}{\wp(h(z))}].\\
\end{aligned}
\end{equation}
Again by taking the residue of both sides of (\ref{X2.12}) at $b_0$, we have
$$
s[\frac{1}{i}+\frac{1}{\sqrt{3}}]=A\eta\frac{m}{\sqrt{3}},
$$
which implies that $A\eta=\frac{s}{m}(1-i\sqrt{3})$. Combining $A\eta=\frac{\sqrt{3}+i}{(\sqrt{3}-i)}$ yields that
$$
\frac{s}{m}=\frac{\sqrt{3}+i}{(\sqrt{3}-i)(1-i\sqrt{3})}=\frac{-1+\sqrt{3}i}{4},
$$
which is a contradiction. Therefore, $\wp'(h(\overline{b_0}))=-i$. Thus, all the poles of $\wp(h(z))$ must be the zero of $\wp'(h(\overline{z}))+i$. Suppose that $\wp(h(z_0))=0$ and $\wp'(h(z_0))=i$. The same argument also yields $\wp(h(\overline{z_0}))=0$ and $\wp'(h(\overline{z_0}))=i$, which means that all the zeros $\wp'(h(z))-i$ must be the zero of $\wp'(h(\overline{z}))-i$. We denote the facts as follows.
\begin{equation}\label{XX1}
\wp(h(z))=\infty\Rightarrow \wp'(h(\overline{z}))+i=0,~~\wp'(h(z))-i=0\Rightarrow \wp'(h(\overline{z}))-i=0.
\end{equation}

Suppose that 0 is a Picard value of $h(z)$. Then, we can assume that $h(z)=e^{\alpha(z)}$, where $\alpha(z)$ is an entire function. The equation $h(\overline{z})=Ah(z)+B$ yields $e^{\alpha(\overline{z})}=Ae^{\alpha(z)}+B=A[e^{\alpha(z)}-(-\frac{B}{A})]$. So, $-\frac{B}{A}$ is also a Picard value of $h(z)=e^{\alpha(z)}$, and Picard's little theorem tells us $B=0$. Thus, $h(\overline{z})=Ah(z)$. \\

Suppose that $\tau(\neq 0)\in \mathbb{L}$, which means that $\tau$ is any fixed period of $\wp$. Note that $\tau$ is not Picard value of $h(z)$. Assume that $h(u_0)=\tau$. Then, $\wp(h(u_0))=\wp(\tau)=\infty$. From (\ref{XX1}), one has $\wp(h(\overline{u_0}))=0$ and $\wp'(h(\overline{u_0}))=-i$. Without loss of generality, we assume that
\begin{equation}\label{XX2}
h(\overline{u_0})=Ah(u_0)=A\tau=\theta_1+\xi_1,
\end{equation}
where $\xi_1\in\mathbb{L}$. We also know that $\theta_2$ is not Picard value of $h(z)$. Assume $h(e_0)=\theta_2$. Thus, $\wp(h(e_0))=\wp(\theta_2)=0$ and $\wp'(h(e_0))=\wp'(\theta_2)=i$. Then, (\ref{XX1}) yields $\wp(h(\overline{e_0}))=0$ and $\wp'(h(\overline{e_0}))=i$.  Without loss of generality, we assume that
\begin{equation}\label{XX3}
h(\overline{e_0})=Ah(e_0)=A\theta_2=\theta_2+\xi_2,
\end{equation}
where $\xi_2\in\mathbb{L}$. Clearly, $\theta_2+\tau$ is also not Picard value of $h(z)$. Assume $h(t_0)=\theta_2+\tau$. Then, $\wp(h(t_0))=\wp(\theta_2+\tau)=0$ and $\wp'(h(t_0))=\wp'(\theta_2+\tau)=i$. The same argument as above yields that
\begin{equation}\label{XX4}
h(\overline{t_0})=Ah(t_0)=A(\theta_2+\tau)=\theta_2+\xi_3,
\end{equation}
where $\xi_3\in\mathbb{L}$. Combining (\ref{XX3}) and (\ref{XX4}) yields that $A\tau=\xi_3-\xi_2\in \mathbb{L}$. Together with (\ref{XX2}) yields that $\theta_1\in \mathbb{L}$ is a period of $\wp$, a contradiction. Thus, $0$ is not a Picard value of $h(z)$. \\

Below, for any period $\tau$ of $\wp$, we will prove that $A\tau\in \mathbb{L}$, which means that $A\tau$ is also a period of $\wp$.\\

Note that $0$ is not a Picard value of $h(z)$. Assume that $h(d_0)=0$. Thus, $\wp(h(d_0))=\wp(0)=\infty$. Then, (\ref{XX1}) yields $\wp(h(\overline{d_0}))=0$ and $\wp'(h(\overline{d_0}))=-i$. Without loss of generality, we assume that
\begin{equation}\label{XX5}
h(\overline{d_0})=Ah(d_0)+B=B=\theta_1+\xi_5.
\end{equation}
where $\xi_5\in\mathbb{L}$. Observe that $h(z)$ has one finite Picard value at most. So, one of $\omega_1$ and $\omega_2$ is not a Picard value of $h(z)$. Without loss of generality, we assume $\omega_1$ is not a Picard value of $h(z)$ and $h(c_0)=\omega_1$. Then, $\wp(h(c_0))=\wp(\omega_1)=\wp(0)=\infty$, which plus (\ref{XX1}) implies that $\wp(h(\overline{c_0}))=0$ and $\wp'(h(\overline{c_0}))=-i$. So, we can set $h(\overline{c_0})=\theta_1+\xi_6$ with $\xi_6\in\mathbb{L}$. Further,
\begin{equation}\label{XX6}
h(\overline{c_0})=Ah(c_0)+B=Ah(c_0)+\theta_1+\xi_5=A\omega_1+\theta_1+\xi_5= \theta_1+\xi_6,
\end{equation}
which implies that $A\omega_1=\xi_6-\xi_5\in\mathbb{L}$.\\

Clearly, there exists an integer $N$ such that $N\omega_1+\tau$ is not a Picard value of $h(z)$. Suppose that $h(t_1)=N\omega_1+\tau$. Then, $\wp(h(t_1))=\wp(N\omega_1+\tau)=\wp(0)=\infty$, which also implies that $\wp(h(\overline{t_1}))=0$ and $\wp'(h(\overline{t_1}))=-i$. We can set $h(\overline{t_1})=\theta_1+\xi_7$ with $\xi_7\in\mathbb{L}$. Further,
\begin{equation}\label{XX7}
h(\overline{c_1})=Ah(c_1)+\theta_1+\xi_5=A(N\omega_1+\tau)+\theta_1+\xi_5 =\theta_1+\xi_7,
\end{equation}
which leads to
$$
A\tau=\xi_7 -\xi_5-NA\omega_1\in\mathbb{L}.
$$
Thus, we prove that $A\tau$ is also a period of $\wp$. \\

Observe that $A=\frac{(1+i\frac{\sqrt{3}}{3})}{\eta(1-i\frac{\sqrt{3}}{3})}$ and $\eta^3=1$. A calculation yields $A^3=-1$ and
$$
A=\frac{1}{2}+\frac{\sqrt{3}}{2}i=e^{\frac{\pi}{3}i},~~or~~\frac{1}{2}-\frac{\sqrt{3}}{2}i=e^{-\frac{\pi}{3}i},~~or~~-1.
$$
Meanwhile, (\ref{XX5}) yields that
$$
h(\overline{z})=Ah(z)+B=Ah(z)+\theta_1+\xi_5,
$$
which is (1) with $\tau_1=\xi_5$.\\

Next, we consider $g_3(z)\equiv 0$ or $g_2(z)\equiv 0$.\\

If $g_3(z)\equiv 0$, as above discussion, then, one can derive
$$
\wp(h(z))=\infty\Rightarrow \wp'(h(\overline{z}))-i=0,~~\wp'(h(z))+i=0\Rightarrow \wp'(h(\overline{z}))+i=0.
$$

If $g_2(z)\equiv 0$, as above discussion, then, one can derive
$$
\wp(h(z))=\infty\Rightarrow \wp(h(\overline{z}))=\infty,~~\wp'(h(z))-i=0\Rightarrow \wp'(h(\overline{z}))+i=0.
$$
Further, with the same argument, we can get the conclusions (2) and (3) if $g_3(z)\equiv 0$ and $g_2(z)\equiv 0$, respectively. Here, we omit the details.\\

Next, we will prove that $L(z)$ is a linear function with two constant $a(\neq0),~b$. Without loss of generality, we assume that
\begin{equation}\label{Z11}
h(L(z))=Ah(z)+C,
\end{equation}
where $A^3=-1$ and $C$ is a constant. Suppose that $h$ is transcendental. We recall the following result (see \cite[Theorem 2]{Clu} and \cite[p. 370]{GY}).\\

\emph{Lemma 2. If $f$ (meromorphic) and $g$ (entire) are transcendental, then
$$
\limsup_{r\not\in E,~r\rightarrow\infty}\frac{T (r, f (g))}{T (r, f )}=\infty,
$$
where $E$ is a set of finite Lebesgue measure.}\\

If $L$ is transcendental, by Lemma 2, we have that
$$
\infty=\limsup_{r\not\in E,~r\rightarrow\infty}\frac{T(r,h(L(z)))}{T(r,h(z))}=\limsup_{r\not\in E,~r\rightarrow\infty}\frac{T(r,Ah(z)+C)}{T(r,h(z))}=1,
$$
a contradiction. Thus, $L$ is a polynomial. Without loss of generality, we assume $L(z)=a_{m}z^{m}+\cdots+a_{1}z+a_{0}$ with $m\geq 1$ and $a_m\neq 0$. We need the following result, which can be seen in \cite[(19)]{GY}) and \cite[(2.7)]{Li3}), respectively. \\

\emph{Lemma 3. Suppose that $g(z) = a_m z^m+ a_{m_1}z^{m-1}+ ... + a_1 z + a_0$ ($a_m\neq 0$) is a non-constant polynomial,
then for any $\epsilon_1 > 0$ and $\epsilon_2 > 0$,
$$
 T (r, f (g))\geq (1 -\epsilon_2 )T(\frac{a_m}{2}r^m, f),
$$
$$
 T (r, f)\leq \frac{1}{m}(1+\epsilon_1)T(\frac{a_m}{2}r^m, f),
$$
for large $r$ outside possibly a set of finite Lebesgue measure. }\\

Applying Lemma 3 to the function $h(L(z))$, we obtain
\begin{equation}\label{Z12}
\begin{aligned}
T(r,h(z))&\leq \frac{1}{m}(1+\epsilon_1)T(\frac{a_m}{2}r^m, h(z))\leq \frac{1}{m}\frac{1+\epsilon_1}{1 -\epsilon_2 }T (r, h(L(z)))\\
&=  \frac{1}{m}\frac{1+\epsilon_1}{1 -\epsilon_2 }T (r, Ah(z)+C))= \frac{1}{m}\frac{1+\epsilon_1}{1 -\epsilon_2 }T (r, h(z))+O(1),
\end{aligned}
\end{equation}
which implies $m=1$, since $\epsilon_1$ and $\epsilon_2$ can be chosen small enough. Thus, $L(z)$ is a linear function. \\

Suppose that $h$ is a polynomial with degree $n$. Then, (\ref{Z11}) yields that $h(L(z))$ is also a polynomial, which implies that $L(z)$ is also a polynomial. Further, (\ref{Z11}) leads to that $L(z)$ is a linear function.\\

Thus, the above discussion yields that $L(z)=qz+c$ with $q\neq0$ and $c$ is a constant. So, $h(qz+c)=Ah(z)+C$, where $C$ is a constant.\\

At the end, we will prove $|q|=1$ by an elementary method. Differentiating the above equation leads to
\begin{equation}\label{N11}
qh'(qz+c)=Ah'(z).
\end{equation}
If $h'(z)$ is a constant, then, $q=A$ and $|q|=|A|=1$. Next, we assume that $h'(z)$ is nonconstant. Suppose that $|q|<1$. For each $t^*\in \mathbb{C}$, by (\ref{N11}), we get
$$
h'(t^*)=\frac{q}{A}h'(qt^*+c)=(\frac{q}{A})^2h'[q(qt^*+c)+c]=\cdots=(\frac{q}{A})^nh'[q^n t^*+c\sum_{k=0}^{n-1}q^{k}].
$$
Note that $|q|<1$ and $q^n t^*+c\sum_{k=0}^{n-1}q^{k}\rightarrow c\frac{1}{1-q}$ as $n\rightarrow\infty$. Then, let $n\rightarrow\infty$, one has that
$$
h'(t^*)=h'(c\frac{1}{1-q})\lim_{n\rightarrow\infty}(\frac{q}{A})^n=0,
$$
which means that $h'(z)\equiv 0$, a contradiction. If $|q|> 1$, then we rewrite (\ref{N11}) as $h'(z)=Aph'(pz+c')$ with $p=\frac{1}{q}$ and $c'=-\frac{c}{q}$. Note that $|p|=|\frac{1}{q}|<1$. The same argument as above yields a contradiction. So, we obtain $|q|=1$.\\

Thus, we finish the proof of this result.

\end{proof}

\begin{proof}[Proof of Proposition 2] The necessity is obvious. Below, we prove the sufficiency. We rewrite (\ref{L}) as $$[\frac{f(z)}{e^{\frac{g(z)}{n}}}]^n+[\frac{f(L(z))}{e^{\frac{g(z)}{m}}}]^m=1.$$
We will prove that both $\frac{f(z)}{e^{\frac{g(z)}{n}}}$ and $\frac{f(L(z))}{e^{\frac{g(z)}{m}}}$ are constant.

Observe that Yang's Theorem yields that both $\frac{f(z)}{e^{\frac{g(z)}{n}}}$ and $\frac{f(L(z))}{e^{\frac{g(z)}{m}}}$ are constant if $\frac{2}{m}+\frac{1}{n}<1$. So, it is suffice to consider the case $(n,m)=(2,3)$ and $(2,4)$. \\

\textbf{Case 1.} $(n,m)=(2,3)$. \\

Suppose that one of $\frac{f(z)}{e^{\frac{g(z)}{2}}}$ and $\frac{f(L(z))}{e^{\frac{g(z)}{3}}}$ is not constant. Then neither of them is constant. Via (ii) of Theorem A, we have that $\frac{f(z)}{e^{\frac{g(z)}{2}}} = i\wp'(h(z))$ and $\frac{f(L(z))}{e^{\frac{g(z)}{3}}} = \eta\sqrt[3]{4}\wp(h(z))$, where $h(z)$ is any nonconstant entire function. So,
\begin{equation}\label{Y111111}
f(L(z)) =i\wp'(h(L(z)))e^{\frac{g(L(z))}{2}}= \eta\sqrt[3]{4}\wp(h(z))e^{\frac{g(z)}{3}}.
\end{equation}
Note that $w_n=n\omega_1$ $(n=0, ~1,~2,...)$ is a pole of $\wp(z)$ with multiplicity 2. It is known that a nonconstant meromorphic function has four complete multiple values at most. (Here, the constant $a$ is a complete multiple value of $f$ if $f-a$ only has multiple zeros). So, there exists a point $w_N$ such that $h(z)-w_N$ has simple zero, say $z_0$. Then, $h(z_0)=w_N$ and $z_0$ is a pole of the function
$\eta\sqrt[3]{4}\wp(h(z))e^{\frac{g(z)}{3}}$ with multiplicity 2. On the other hand, in view of that all the poles of $\wp'(w)$ is 3, we derive that $z_0$ is a pole of $i\wp'(h(L(z)))e^{\frac{g(L(z))}{2}}$ at least 3 or $z_0$ is not a pole of $i\wp'(h(L(z)))e^{\frac{g(L(z))}{2}}$.  Then comparing the multiplicities of both sides of (\ref{Y111111}) at pole-point $z_0$, we have a contradiction.\\

\textbf{Case 2.} $(n,m)=(2,4)$. \\

we assume that one of $\frac{f(z)}{e^{\frac{g(z)}{2}}}$ and $\frac{f(L(z))}{e^{\frac{g(z)}{4}}}$ is not a constant. Via (iii) of Theorem A, we have that $\frac{f(z)}{e^{\frac{g(z)}{2}}} = sn'(h(z))$ and $\frac{f(L(z))}{e^{\frac{g(z)}{4}}} = sn(h(z))$, where $h(z)$ is any nonconstant entire function and $sn$ is Jacobi elliptic function satisfying ${sn'}^2=1- sn^4$. So,
\begin{equation}\label{Y}
f(L(z)) =sn'(h(L(z)))e^{\frac{g(L(z))}{2}}= sn(h(z))e^{\frac{g(z)}{4}}.
\end{equation}
Suppose that $w_n$ $(n=0, ~1,~2,...)$ is a simple pole of $sn(z)$. Similarly as above case, there exists a point $w_N$ such that $h(z)-w_N$ has simple zero, say $z_1$. Then, $h(z_1)=w_N$ and $z_1$ is a pole of the function $sn(h(z))e^{\frac{g(z)}{4}}$ with multiplicity 1. On the other hand, in view of that all the pole of $sn'(w)$ is 2, we get that $z_1$ is a pole of $sn'(h(z))e^{\frac{g(L(z))}{2}}$ at least 2 or $z_1$ is not a pole of $i\wp'(h(L(z)))e^{\frac{g(L(z))}{2}}$. Then comparing the multiplicities of both sides of (\ref{Y}) at pole-point $z_1$, we have a contradiction.\\

Therefore, we derive that both $\frac{f(z)}{e^\frac{g(z)}{n}}$ and $\frac{f(L(z))}{e^\frac{g(z)}{m}}$ are nonzero constants, say $A$ and $B$ with $A^n+B^m=1$. Thus, $f(z)=Ae^\frac{g(z)}{n}$ and
$$
f(L(z))=Be^\frac{g(z)}{m}=Ae^\frac{g(L(z))}{n},
$$
which implies that
\begin{equation}\label{K3}
g(L(z))=\frac{n}{m}g(z)+nLn\frac{B}{A}.
\end{equation}
With the same argument as in Proposition 1, we can derive that $L(z)$ is a linear function, say $L(z)=qz+c$ with constants $q(\neq0)$, $c$. \\

If $m=n$, then (\ref{K3}) reduces to $g(qz+c)=g(z)+nLn\frac{B}{A}$. Then, the same discussion as in Proposition 1 yields $|q|=1$. \\

Thus, we finish the proof of this result.

\end{proof}

Based on Theorem 1, we give the proof of Theorem 2 as follows.\\

\begin{proof}[Proof of Theorem 2] The necessity is obvious. Below, we prove the sufficiency. Suppose that the hyper-order of $f(z)$ is less than 1. Then by Theorem B, we get the conclusion (a). Next, we assume that the hyper-order $\rho_{2}(f)\geq 1$. Rewrite (\ref{Y11}) as
$$
(\frac{f(z+c)}{e^{\frac{P(z)}{3}}})^3+(\frac{f(z)}{e^{\frac{P(z)}{3}}})^3=1.
$$
Via (i) of Theorem A, we have
\begin{equation}\label{Y2}
f(z)=\frac{1}{2}\frac{1+\frac{\wp'(h(z))}{\sqrt{3}}}{\wp(h(z))}e^{P(z)/3}
\end{equation}
where $h(z)$ is an entire function. The fact $\rho_{2}(f)\geq 1$ and (\ref{Y2}) yields that $h$ is transcendental. Here, we employ a result of Bergweiler in \cite[Lemma 1]{Ber}.\\

\emph{Lemma 3. Let $f$ be meromorphic and let $g$ be entire and transcendental. If the lower order $\mu(f \circ g)<\infty,$ then $\mu(f)=0$.}\\

Note that $\rho(\wp)=\mu(\wp)=2$, see \cite{B.L}. It follows from Lemma 3 that $\mu[\frac{1}{2}\frac{1+\frac{\wp'(h(z))}{\sqrt{3}}}{\wp(h(z))}]=\infty$, so is $\mu(f)$. Thus $e^P$ is a small function of $f(z)$ and $T(r,e^P)=S(r,f)$. Based on the idea of Liu and Ma in \cite{LM}, we will obtain the desired result.\\

Set $a(z)=e^{P(z)}$ and rewrite (\ref{Y11}) as $f^{3}(z+c)=-(f^{3}(z)-a(z))$, which implies that the zeros of $f^{3}(z)-a(z)$ are of multiplicities at least 3. Rewrite (\ref{Y11}) as $f^{3}(z)=-(f^{3}(z+c)-a(z))$, which implies that the zeros of $f^{3}(z+c)-a(z)$ are of multiplicities at least 3. So, the zeros of $f^{3}(z)-a(z-c)$ are of multiplicities at least 3. Set $G(z)=f^{3}(z)$. Assume that the functions $a(z-c)$ and $a(z)$ are distinct from each other. Then, applying the second main theorem of Nevanlinna to $G$, one gets that
$$
\begin{aligned}
2T(r,G)&\leq \overline{N}(r,\frac{1}{G(z)-a(z-c)})+\overline{N}(r,\frac{1}{G(z)-a(z)})+\overline{N}(r,G)+\overline{N}(r,\frac{1}{G})\\
&+S(r,G)\leq\frac{1}{3}[N(r,\frac{1}{G(z)-a(z-c)})+N(r,\frac{1}{G(z)-a(z)})+N(r,G)\\
&+N(r,\frac{1}{G})]+S(r,G)\leq \frac{4}{3}T(r,G)+S(r,G),
\end{aligned}
$$
which is a contradiction. Thus, $a(z)=a(z-c)$, which implies that $e^{P(z+c)}=e^{P(z)}$ and $e^{P(z+c)-P(z)}=1$. In view of that $P(z)$ is a polynomial, we derive that $P(z)=\alpha z+\beta$ with two constants $\alpha, ~\beta$ and $e^{\alpha c}=1$. Therefore, we can rewrite (\ref{Y11}) as
$$
1=[\frac{f(z)}{e^{\frac{P(z)}{3}}}]^3+[\frac{f(z+c)}{e^{\frac{P(z)}{3}}}]^3=[\frac{f(z)}{e^{\frac{P(z)}{3}}}]^3+[\frac{f(z)}{e^{\frac{P(z+c)}{3}}}]^3.
$$
Set $F(z)=\frac{f(z)}{e^{\frac{P(z)}{3}}}$. Then, $F(z)^3+F(z+c)^3=1$. By Theorem 1, we get the desired result.\\

Thus, we finish the proof of Theorem 2.
\end{proof}

At the end of this section, we give the proof of Theorem 3. \\

\begin{proof}[Proof of Theorem 3] It is suffice to consider the case $\rho(g(z))<1$. Suppose that $f$ is a nonconstant meromorphic solution of (\ref{L1}), then $g$ is nonconstant, since $f(z)=Ae^\frac{g(z)}{n}$. Differentiating the equation $g(z+c)=\frac{n}{m}g(z)+nLn\frac{B}{A}$ one time yields that
\begin{equation}\label{K4}
g'(z+c)=\frac{n}{m}g'(z).
\end{equation}
Assume that $g'$ is not a constant. If 0 is a Picard value of $g'(z)$, then, $g'(z)=e^{h(z)}$, where $h(z)$ is a nonconstant entire function. So, the order $\rho(g(z))=\rho(g'(z))=\rho(e^{h(z)})\geq 1$, a contradiction. Therefore, there exists a point $z_2$ such that $g'(z_2)=0$. The above equation yields that $z_2+nc$ also zero of $g'$. So, $n(d|c|+|z_2|,g')\geq d$ for any $d\in \mathbb{N}$, where $n(r,g')$ denotes the number of zeros of $g'$ in $\{z:|z|<r\}$. Based on the method in \cite[Lemma 3.2]{HAL1}, we will derive a contradiction as follows.
$$
\begin{aligned}
\rho(g(z))&=\rho(g'(z))\geq \limsup _{r \rightarrow \infty} \frac{\log n(r, g')}{\log r}\geq \limsup _{d \rightarrow \infty} \frac{\log n(d|c|+|z_2|, g')}{\log (d|c|+|z_2|)}\\
&\geq \limsup _{d \rightarrow \infty} \frac{\log d}{\log (d|c|+|z_2|)}=1 ,
\end{aligned}
$$
a contradiction. Thus, $g'$ is a nonzero constant and (\ref{K4}) yields that $m=n$. Further, $g(z)=\alpha z+\beta$, and $f(z)=Ae^\frac{\alpha z+\beta}{n}$ with $A^n(1+e^{\alpha c})=1$. \\

Thus, we finish the proof of Theorem 3.

\end{proof}


\begin{thebibliography}{99}

\bibitem{Ba}I.N. Baker, On a class of meromorphic functions, Proc. Amer. Math. Soc. 17(1966), 819-822.

\bibitem{B.L}S.B. Bank and J.K. Langley, On the value distribution theory of elliptic functions, Monatsh. Math. 98(1984), 1-20.

\bibitem{Bar}D.C. Barnett, R.G. Halburd, R.J. Korhonen and W. Morgan, Nevanlinna theory for the $q$-difference operator and meromorphic solutions of $q$-difference equations, Proc. R. Edinburgh. A. 137(2007), 457-474.

\bibitem{Ber}W. Bergweiler, Order and lower order of composite meromorphic functions, Mich. Math. J. 36(1989), 135-146.

\bibitem{Ber1}W. Bergweiler and W. Hayman, Zeros of solutions of a functional equation, Comput.
Methods Funct. Theory. 3(2003), 55-78.

\bibitem{Ber2}W. Bergweiler, K. Ishizaki, and N. Yanagihara, Meromorphic solutions of some functional
equations, Methods Appl. Anal. 5(1998), 248-258, Correction: Methods Appl. Anal. 6(1999), 617-618.

\bibitem{BL}W.Q. Bi and F. L\"{u}, On meromorphic solutions of the Fermat-type functional equations $f^3(z)+f^3(z+c)=e^{P}$, Anal. Math. Phys. 13(2023), Paper No. 24, 15 pp.

\bibitem{Chen}W. Chen, Q. Han and J.B. Liu, On Fermat Diophantine functional equations, little Picard
theorem and beyond, Aequat. Math. 93(2019), 425-432.

\bibitem{CF}Y.M. Chiang and S.J. Feng, On the Nevanlinna characteristic of $f(z+\eta)$ and difference
equations in the complex plane, Ramanujan J. 16(2008), 105-129.

\bibitem{Clu1}J. Clunie, On integral and meromorphic functions, J. London Math. Soc. 37(1962), 17-27.

\bibitem{Clu}J. Clunie, The composition of entire and meromorphic functions, Mathematical Essays Dedicated to A.J. Macintyre,
OH. Univ. Press, Athens, Ohio, 1970.

\bibitem{Di}L.E. Dickson, History of the theory of numbers, Volume II, Washington Carnegie Institution of Washington, Washington, 1920.

\bibitem{Ere}A. E. Eremenko and M. L. Sodin, Iterations of rational functions and the distribution of
the values of Poincar\'{e} functions, Teor. Funktsii Funktsional. Anal. i Prilozhen. 53(1990), 18-25,
(Russian); translation in J. Soviet Math. 58(1992), 504-509.

\bibitem{G}M.L. Green, Some Picard theorems for holomorphic maps to algebraic varieties, Am. J. Math. 97(1975),
43-75.

\bibitem{FG}F. Gross, On the equation $f^{n}+g^{n} = 1$ II, Bull. Amer. Math. Soc. 74(1968), 647-648.

\bibitem{F.G}F. Gross, On the equation $f^{n} + g^{n} = h^{n}$, Amer. Math. Monthy. 73(1966), 1093-1096.

\bibitem{GY}F. Gross and C. C. Yang, On meromorphic solution of a certain class of
functional-differential equations, Annales Polonici Mathematici. 27(1973), 305-311.

\bibitem{Gun}G. G. Gundersen, J. Heittokangas, I. Laine, J. Rieppo, and D. Yang, Meromorphic solutions
of generalized Schr\"{¨o}der equations, Aequationes Math. 63(2002), 110-135.

\bibitem{HK}R.G. Halburd and R.J. Korhonen, Difference analogue of the lemma on the logarithmic derivative with applications to difference equations, J. Math. Anal. Appl. 314(2006), 477-487.

\bibitem{HAL1}R.G. Halburd and R.J. Korhonen, Growth of meromorphic solutions of delay differential equations, Proc. Amer. Math. Soc. 145(2017), 2513-2526.

\bibitem{HAL}R.G. Halburd, R.J. Korhonen and K. Tohge, Holomorphic curves with shift-invariant hyperplane preimages, Trans.Amer. Math. Soc. 366(2014), 4267-4298.

\bibitem{LH}Q. Han and F. L\"{u}, On the equation $f^{n}(z)+g^{n}(z)=e^{\alpha z+\beta}$, J. Contemp. Mathemat. Anal. 54(2019), 98-102.

\bibitem{Hay}W.K. Hayman, Meromorphic Functions, Clarendon Press, Oxford, 1964.

\bibitem{Is}K. Ishizaki and N. Yanagihara, Borel and Julia directions of meromorphic Schr\"{o}der functions, Math. Proc. Camb. Phil. Soc. 139(2005), 139-147.

\bibitem{Ja}A.V. Jategaonkar, Elementary proof of a theorem of P. Montel on entire functions, J. Lond. Math. Soc.
40(1965), 166-170.

\bibitem{KW}R. Korhonen and Z.T. Wen, Existence of zero-order meromorphic solutions in detecting q-difference Painlev\'{e} equations, Trans. Amer. Math. Soc. 368(2016), 4993-5008.

\bibitem{KZ} R. Korhonen and Y.Y. Zhang, Existence of meromorphic solutions of first order difference
equations, Constr. Approx. 51(2020), 465-504.

\bibitem{LI}I. Laine, Nevanlinna theory and complex differential equations, de Gruyter, Berlin, 1993.

\bibitem{Li3}B. Q. Li and E. G. Saleeby, On solutions of functional-differential equations $f'(x)=a(x)f(g(x))+b(x)f(x)+c(x)$ in the large, Israel J. Math. 162(2007), 335-348.

\bibitem{Li1}B.Q. Li, On meromorphic solutions of $f^2 + g^2 = 1$, Math. Z. 258(2008), 763-771.

\bibitem{Li2}B.Q. Li, On entire solutions of Fermat type partial differential equations, Int. J. Math. 15(2004), 473–485.

\bibitem{Liu5}K. Liu and T.B. Cao, Entire solutions of Fermat type difference differential equations, Electron. J.
Diff. Equ. 59(2013), 10pp.

\bibitem{Liu3}K. Liu, T.B. Cao and H.Z. Cao, Entire solutions of Fermat type differential-difference equations, Arch. Math. 99(2012), 147-155.

\bibitem{Liu4}K. Liu and L.Z. Yang, On entire solutions of some differential–difference equations, Comput. Methods Funct. Theory. 13(2013), 433-447.

\bibitem{LM}K. Liu and L. Ma, Fermat type equations or systems with composite functions, Journal of Computational Analysis and Applications. 26(2019), 362-372.

\bibitem{LG}F. L\"{u} and H.X. Guo, On the Fermat-type equation $f^{3}(z)+f^{3}(z+c)=e^{\alpha z+\beta }$, Mediterr. J. Math. 19(2022), Paper No. 118, 13 pp.

\bibitem{FL}F. L\"{u} and Q. Han, On the Fermat-type equation $f^{3}(z)+f^{3}(z+c)=1$, Aequat. Math. 91(2017), 129-136.

\bibitem{PM}P. Montel, Le.cons sur les familles normales de fonctions analytiques et leurs applications, G authier-Villars, Paris. 32(1927), 135-136.

\bibitem{O}M. Ozawa, On the existence of prime periodic entire functions, Kodai. Math. Sem. Rep. 29(1978), 308-321.

\bibitem{W2}R. Taylor and A. Wiles, Ring-theoretic properties of certain Hecke algebras, Ann. Math. 141(1995), 553-572.

\bibitem{V}G. Valiron, Fonctions Analytiques, Press. Univ. de France, Paris, 1952.

\bibitem{W1}A. Wiles, Modular elliptic curves and Fermat's last theorem, Ann. Math. 141(1995), 443-551.

\bibitem{WL}L. Wu, C. He, W.R. L\"{u} and F. L\"{u}, Existence of meromorphic solutions of some generalized Fermat functional equations, Aequat. Math. 94(2020), 59-69.

\bibitem{Yana}N. Yanagihara, Polynomial difference equations which have meromorphic solutions of finite order, Analytic function theory of one complex variable, Pitman Res. Notes Math. Ser. 212(1989), 368-392.

 \bibitem{Y}C.C. Yang, A generalization of a theorem of P. Montel on entire functions, Proc. Am. Math. Soc. 26(1970), 332-334.

\bibitem{YY}C.C. Yang and H.X. Yi, Uniqueness Theory of Meromorphic Functions, Science Press, Beijing/New York, 2003.

















\end{thebibliography}
\end{document}